\documentclass[10pt]{article}

\usepackage{array}
\usepackage{multirow}

\makeatletter
\newcommand{\thickhline}{%
    \noalign {\ifnum 0=`}\fi \hrule height 1pt
    \futurelet \reserved@a \@xhline
}
\newcolumntype{"}{@{\hskip\tabcolsep\vrule width 1pt\hskip\tabcolsep}}
\makeatother

\usepackage{graphics}
\usepackage{epsfig}
\usepackage{subfigure}
\usepackage{amssymb}
\usepackage{amsthm}
\usepackage{mathrsfs}
\usepackage{hhline}
\usepackage{longtable}
\usepackage{amsmath}
\usepackage{array}
\usepackage{float}
\usepackage{picinpar}
\usepackage{multirow}
\usepackage{cite}
\usepackage{enumerate}
\usepackage{xcolor}
\usepackage[mathlines]{lineno}
\modulolinenumbers[2]
\newcommand*\patchAmsMathEnvironmentForLineno[1]{%
\expandafter\let\csname old#1\expandafter\endcsname\csname #1\endcsname  \expandafter\let\csname oldend#1\expandafter\endcsname\csname end#1\endcsname  \renewenvironment{#1}%
{\linenomath\csname old#1\endcsname}%
{\csname oldend#1\endcsname\endlinenomath}}%
\newcommand*\patchBothAmsMathEnvironmentsForLineno[1]{%
\patchAmsMathEnvironmentForLineno{#1}%
\patchAmsMathEnvironmentForLineno{#1*}}%
\AtBeginDocument{%
\patchBothAmsMathEnvironmentsForLineno{equation}%
\patchBothAmsMathEnvironmentsForLineno{align}%
\patchBothAmsMathEnvironmentsForLineno{flalign}%
\patchBothAmsMathEnvironmentsForLineno{alignat}%
\patchBothAmsMathEnvironmentsForLineno{gather}%
\patchBothAmsMathEnvironmentsForLineno{multline}%
}

\def\nt{\noindent}
\def\ms{\medskip}

\newtheorem{theorem}{Theorem}[section]
\newtheorem{lemma}[theorem]{Lemma}
\newtheorem{example}{Example}[section]
\newtheorem{corollary}[theorem]{Corollary}
\newtheorem{problem}{Problem}[section]

\newtheorem{conjecture}{Conjecture}[section]

\numberwithin{equation}{section}
\def\Z{\mathbb{Z}}

\newtheorem{defi}{Definition}[section]

\usepackage{graphicx}
\graphicspath{%
    {converted_graphics/}
    {/}
}
\begin{document}
\baselineskip18truept
\normalsize
\begin{center}
{\mathversion{bold}\Large \bf On local antimagic chromatic number of cycle-related join graphs II}

\bigskip
{\large  Gee-Choon Lau{$^{a,}$}\footnote{Corresponding author.}, K. Premalatha{$^b,$} S. Arumugam{$^b,$} Wai-Chee Shiu{$^c$}}\\

\medskip

\emph{{$^a$}Faculty of Computer \& Mathematical Sciences,}\\
\emph{Universiti Teknologi MARA, Johor Branch, Segamat Campus,}\\
\emph{85000 Malaysia.}\\
\emph{geeclau@yahoo.com}\\

\medskip

\emph{{$^b$}National Centre for Advanced Research in Discrete Mathematics,}\\
\emph{Kalasalingam Academy of Research and Education,}\\
\emph{ Anand Nagar, Krishnankoil-626 126, Tamil Nadu, India.}\\
\emph{premalatha.sep26@gmail.com}\\

\medskip

\emph{{$^c$}Department of Mathematics, The Chinese University of Hong Kong,}\\
\emph{Shatin, Hong Kong, P.R. China.}\\
\emph{wcshiu@hkbu.edu.hk}\\

\end{center}


\title{On local antimagic (total) chromatic number of join graphs}

\medskip


\begin{abstract}
An edge labeling of a graph $G = (V, E)$ is said to be local antimagic if it is a bijection $f:E \to\{1,\ldots ,|E|\}$ such that for any pair of adjacent vertices $x$ and $y$, $f^+(x)\not= f^+(y)$, where the induced vertex label of $x$ is $f^+(x)= \sum_{e\in E(x)} f(e)$ ($E(x)$ is the set of edges incident to $x$).  The local antimagic chromatic number of $G$, denoted by $\chi_{la}(G)$, is the minimum number of distinct induced vertex labels over all local antimagic labelings of $G$. In this paper, several sufficient conditions to determine the local antimagic chromatic number of the join of graphs are obtained. We then determine the exact value of the local antimagic chromatic number of many join graphs.
\end{abstract}

\section{Introduction}
For a connected graph $G = (V, E)$, let $f: E \to \{1,\ldots, |E|\}$ be a bijection such that for each vertex $u\in V$, the weight $w(u) = \sum_{e\in E(u)} f(e)$, where $E(u)$ is the set of edges incident to $u$. We say $G$ is antimagic if $w(u)\ne w(v)$ for every two distinct vertices $u,v\in V$, and $f$ is called an antimagic labeling of $G$. Many papers have been published (see, for examples~\cite{Alon, Chang+L+P+Z, Cranston, Cranston+L, Eccles}) since the introduction of the concept of antimagic labeling by Hartsfield~\cite{Hartsfield}. Two of the most famous open conjectures raised by Hartsfield are:

\begin{conjecture} Every connected graph other than $K_2$ is antimagic. \end{conjecture}

\begin{conjecture} Every tree other than $K_2$ is antimagic. \end{conjecture}

\nt A connected graph $G = (V, E)$ is said to be {\it local antimagic} if it admits a {\it local antimagic edge labeling}, i.e., a bijection $f : E \rightarrow \{1,\dots ,|E|\}$ such that the induced vertex labeling $f^+ : V \rightarrow \Z$ given by $f^+(u) = \sum_{e\in E(u)} f(e)$ has the property that any two adjacent vertices have distinct induced vertex labels (see~\cite{Arumugam, Bensmail}). Thus, $f^+$ is a coloring of $G$. Clearly, the order of $G$ must be at least 3.  The vertex label $f^+(u)$ is called the {\it induced color} of $u$ under $f$ (the {\it color} of $u$, for short, if no ambiguity occurs). The number of distinct induced colors under $f$ is denoted by $c(f)$, and is called the {\it color number} of $f$. The {\it local antimagic chromatic number} of $G$, denoted by $\chi_{la}(G)$, is $\min\{c(f) : f\ \mbox{is a local antimagic labeling of}\ G\}$. It is proved~\cite{Haslegrave} that every connected graph except $K_2$ is local antimagic.

\ms\nt Throughout this paper, let the path of order $m$ be $P_m= u_1u_2\cdots u_m$ and $O_n$ be the null graph of order $n\ge 1$ with vertices $v_j, 1\le j\le n$. For any graph $G$ and $H$, the join graph $G\vee H$ is defined by $V(G\vee H)=V(G)\cup V(H)$ and $E(G\vee H)=E(G) \cup E(H) \cup \{uv \,|\, u \in V(G), v\in V(H)\}$. For integers $a < b$, let $[a,b]=\{a,a+1,\ldots,b\}$. Notations not defined here are referred to the book by Bondy and Murty~\cite{Bondy+M}. 

\ms\nt The following lemmas in~\cite{LSN} are needed. 

\begin{lemma}\label{lem-2part} Let $G$ be a graph of size $q$. Suppose there is a local antimagic labeling of $G$ inducing a $2$-coloring of $G$ with colors $x$ and $y$, where $x<y$. Let $X$ and $Y$ be the number of vertices of color $x$ and color $y$, respectively. Then $G$ is a bipartite graph whose sizes of parts are $X$ and $Y$ with $X>Y$, and
\begin{equation}\label{eq-bi} xX=yY= \frac{q(q+1)}{2}.\end{equation}
\end{lemma}

\ms\nt The contrapositive of Lemma~\ref{lem-2part} gives a sufficient condition for a bipartite graph $G$ to have $\chi_{la}(G)\ge 3$. Let $G-e$ be the graph $G$ with an edge $e$ deleted. 

\begin{lemma}\label{lem-reg} Suppose $G$ is a $d$-regular graph of size $q$. If $f$ is a local antimagic labeling of $G$, then $g = q + 1 - f$ is also a local antimagic labeling of $G$ with $c(f)= c(g)$.  Moreover, if $c(f)=\chi_{la}(G)$ and $f(uv)=1$ or $f(e)=q$, then $\chi_{la}(G-e)\le \chi_{la}(G)$. \end{lemma}



\nt Note that if $G$ is a regular edge-transitive graph, then $\chi_{la}(G-e)\le \chi_{la}(G)$.

\begin{lemma}\label{lem-nonreg} Suppose $G$ is a graph of size $q$ and $f$ is a local antimagic labeling of $G$. For any $x,y\in V(G)$, if\\
(i) $f^+(x) = f^+(y)$ implies that $\deg(x)=\deg(y)$, and\\
(ii) $f^+(x) \ne f^+(y)$ implies that $(q+1)(\deg(x)-\deg(y)) \ne f^+(x) - f^+(y)$,\\ then $g = q + 1 - f$ is also a local antimagic labeling of $G$ with $c(f)= c(g)$. \end{lemma}


\nt Let $G$ be a graph with $\chi_{la}(G)=t\geq 2.$ Let $f$ be a local antimagic labeling with induces a $t$-coloring of $G.$ Let $V_i=\{v\in V:$ The color of $v$ induced by $f$ is $i\}.$ Clearly each $V_i$ is an independent set and $V_f=\{V_1,V_2,\dots,V_t\}$ is a partition of $V(G).$ Suppose the following conditions are satisfied.
\begin{enumerate}[(i)]
\item For each $x\in V_k, 1\leq k\leq t, deg(x)=d_k$ and $f^+(x)-d_a\neq f^+(y)-d_b$ where $x\in V_a, y\in V_b$ and $1\leq a\neq b\leq t.$
\item $f^+(x)+d_a\neq f^+(y)+d_b$ where $x\in V_a, y\in V_b$ and $1\leq a\neq b\leq t.$
\end{enumerate}


\begin{lemma}\label{lem-G-e} Let $H$ be obtained from $G$ with an edge $e$ deleted. If $G$ satisfies Conditions (i) and (ii) and $f(e)=1$, then $\chi(H)\le \chi_{la}(H)\le t$. \end{lemma}



\section{Main Results}

\ms\nt In~\cite[Theorem 3]{LNS}, the authors gave sufficient conditions for $\chi_{la}(G\vee O_n) = \chi_{la}(G)+1$. From the given proof, we can rewrite the theorem as follows.

\begin{theorem}\label{thm-GVOn} Let $G$ be a graph of order $m\geq 3$ and size $e$ with $\chi(G)=\chi_{la}(G)=t.$ Let $f$ be a local antimagic labeling of $G$ that induces a $t$-coloring of $G.$ Let $n\geq 2$ and $m\equiv n(mod\ 2).$ If $f^+(u)\neq (m-n)(2e+mn+1)/2$ for each $u\in V(G),$ then $\chi_{la}(G \vee O_n)=\chi_{la}(G)+1.$
\end{theorem}

\nt By a similar idea, the next two theorems give sufficient conditions for $\chi_{la}(G \vee K_{m,n}) = \chi_{la}(G)+2$ and $\chi_{la}(G \vee C_m) = \chi_{la}(G)+3$.

\begin{theorem}\label{thm-GVKmn} Let $G$ be a graph of even order $p$ and size $e.$ Let $f$ be a local antimagic labeling of $G$ that induces a $\chi(G)$-coloring of $G.$ Suppose $m\neq n\geq 2, m\equiv n(mod\ 2)$ and $(m,n)\neq (2,2).$ If $f^+(u) \not\in \{(p-m)e + (p-m-n)(p(m+n)+1)/2 + np(m+n) + n(mn+1)/2, (p-n)e + (p-m-n)(p(m+n)+1)/2 + mp(m+n) + m(mn+1)/2\}$, then $\chi_{la}(G \vee K_{m,n}) = \chi_{la}(G) + 2$.
\end{theorem} 

\begin{proof} Note that $\chi(G)=\chi_{la}(G)$. Let $V(G) = \{u_i\,|\,1\le i\le p\}$, and $V(K_{m,n}) = \{x_j, y_k\,|\, 1\le j\le m, 1\le k\le n\}$. Suppose $f$ is a local antimagic labeling of $G$ that induces a $\chi(G)$-coloring of $G$. Let $M$ be a magic $(p, m+n)$-rectangle with $(i,j)$-entry $a_{i,j}$, row constant $(m+n)(p(m+n)+1)/2$ and column constant $p(p(m+n)+1)/2$. Also let $N$ be a magic $(m,n)$-rectangle with $(j,k)$-entry $b_{j,k}$, row constant $n(mn+1)/2$ and column constant $m(mn+1)/2$. 

\ms\nt Define $g: E(G\vee K_{m,n})\to [1,e+p(m+n)+mn]$ by $g(uv) = f(uv)$ for each $uv\in E(G)$, 
\begin{eqnarray*}
g(u_ix_j) &=& a_{i,j} + e \mbox{ for } 1\le i\le p, 1\le j\le m, \\
g(u_iy_k) &=& a_{i,m+k} + e \mbox{ for } 1\le k\le n,
\end{eqnarray*}
and  $g(x_jx_k) = e + p(m+n) + b_{j,k}$ for $1\le j\le m, 1\le k\le n$.

\ms\nt It is clear that $g$ is a bijection such that 
\begin{enumerate}[(1)]
  \item $g^+(u) = f^+(u) + (m+n)e + (m+n)(p(m+n)+1)/2$ for each $u\in V(G)$,
  \item $g^+(x_j) = pe + p(p(m+n)+1)/2 + ne + np(m+n) + n(mn+1)/2$,
  \item $g^+(y_k) = pe + p(p(m+n)+1)/2 + me + mp(m+n) + m(mn+1)/2$.
\end{enumerate} 

\nt Since $m\ne n$, we have $(2)\ne (3)$. Since $f^+(u) \not\in\{(p-m)e + (p-m-n)(p(m+n)+1)/2 + np(m+n) + n(mn+1)/2, (p-n)e + (p-m-n)(p(m+n)+1)/2 + mp(m+n) + m(mn+1)/2\}$, we have $g^+(u)\ne g^+(x_j) \ne g^+(x_k)$ for each $u\in V(G)$. Thus, $g$ is a local antimagic labeling that induces $(\chi(G)+2)$-coloring of $G \vee K_{m,n}$. Hence, $\chi_{la}(G\vee K_{m,n}) \le \chi(G) + 2$. Since $\chi_{la}(G\vee K_{m,n}) \ge \chi(G\vee K_{m,n}) = \chi(G)+2$, the theorem holds.
\end{proof}

\begin{theorem}\label{thm-GVC} Let $G$ be a graph of odd order $p$ and size $e.$ Let $f$ be a local antimagic labeling of $G$ that induces a $\chi(G)$-coloring of $G.$ Suppose $m\geq 3$ and $f^+(u) \not\in \{(p-m)(2e+pm+1)/2 + 2(e+pm)+m, (p-m)(2e+pm+1)/2 + 2(e+pm)+m+1, (p-m)(2e+pm+1)/2 + 2(e+pm)+(3m+1)/2\}$, then $\chi_{la}(G \vee C_m) = \chi_{la}(G) + 3$. 
\end{theorem} 

\begin{proof} Note that $\chi(G)=\chi_{la}(G)$. Let $V(G) = \{u_i\,|\,1\le i\le p\}$ and $C_m=v_1v_2\cdots v_mv_1$. Suppose $f$ is a local antimagic labeling of $G$ that induces a $\chi(G)$-coloring of $G$. Let $M$ be a magic $(p, m)$-rectangle with $(i,j)$-entry $a_{i,j}$, row constant $m(pm+1)/2$ and column constant $p(pm+1)/2$.

\ms\nt Define $g : E(G\vee C_m) \to [1,e+pm+m]$ by $g(uv) = f(uv)$ for each $uv\in E(G)$, $g(u_iv_j) = a_{i,j} + e$ for $1\le i\le p, 1\le j\le m$, and $$g(v_iv_{i+1}) = \begin{cases} e+pm+i/2  & \mbox{ if $i$ is even},  \\ e+pm+m-(i-1)/2  & \mbox{ if $i$ is odd (with $v_{p+1}=v_1)$}. \end{cases}$$

\ms\nt It is clear that $g$ is a bijection such that 
\begin{enumerate}[(1)]
  \item $g^+(u) = f^+(u) + me + m(pm+1)/2$ for each $u\in V(G)$,
  \item $g^+(v_i) = pe + p(pm+1)/2 + 2(e+pm) + m$ for odd $i \ge 3$,
  \item $g^+(v_i) = pe + p(pm+1)/2 + 2(e+pm) + m+1$ for even $i$,
  \item $g^+(v_1) = pe + p(pm+1)/2 + 2(e+pm) + (3m+1)/2$.
\end{enumerate}  
Clearly, $(2) < (3) < (4)$. Since $f^+(u) \not\in \{(p-m)(2e+pm+1)/2 + 2(e+pm)+m, (p-m)(2e+pm+1)/2 + 2(e+pm)+m+1, (p-m)(2e+pm+1)/2 + 2(e+pm)+(3m+1)/2\}$, we have $g^+(u)\ne g^+(v)$ for each $u\in V(G)$ and $v\in V(H)$. Thus, $g$ is a local antimagic labeling that induces a $(\chi(G)+3)$-coloring of $G \vee C_m$. Hence, $\chi_{la}(G\vee C_m) \le \chi(G)+3$. Since $\chi_{la}(G\vee C_m) \ge \chi(G\vee C_m) = \chi(G)+3$, the theorem holds.
\end{proof}

\nt In~\cite{LSN}, the authors proved that $\chi_{la}(P_3\vee O_{n})=3, n\ge 3$. They posed the following problem.

\begin{problem}\label{pbm-PmVOn} Determine $\chi_{la}(P_m\vee O_n)$ for $m\ge 4, n\ge 2$. \end{problem}

\nt The following two theorems completely solve Problem~\ref{pbm-PmVOn} when $m$ is even. 

\begin{theorem}\label{thm-P2mVO2n} For $m,n\ge 1$, $\chi_{la}(P_{2m}\vee O_{2n}) = 3$.   \end{theorem}  

\begin{proof} When $m=1$, we get $P_2\vee O_{2n} = K_{1,1,2n}$. In~\cite[Theorem 2.7]{LNS-MJM}, the authors proved that $\chi_{la}(K_{1,1,2n}) = 3$. When $n=1$, the authors in~\cite[Theorem 2.4]{YBYL} proved that $\chi_{la}(P_{2m} \vee O_2) = 3$ for $m\ge 2$. For completeness, the labeling function $f$ is given below:

$$f(u_iu_{i+1}) = \begin{cases} 2m - (i+1)/2 & \mbox{ if $i$ is odd,} \\ i/2 &\mbox{ if $i$ is even.}\end{cases}$$
$$f(u_iv_1) = \begin{cases} 2m + (i-1)/2 & \mbox{ if $i$ is odd,} \\ 6m - (i+2)/2 & \mbox{ if $i$ is even, $i\ne 2m$,} \\ 6m-1 & \mbox{ if $i=2m$.} \end{cases}$$
$$f(uv_2) = \begin{cases} 5m - (i+1)/2 & \mbox{ if $i$ is odd,} \\ 3m + (i-2)/2 & \mbox{ if $i$ is even.} \end{cases}$$

\nt The induced vertex labels are $f^+(u_i) = 9m-2$ for odd $i$, $f^+(u_i) = 11m-2$ for even $i$ and $f^+(v_1) = f^+(v_2) = 8m^2 - m$.

\ms\nt We now consider $m,n\ge 2$. Since $P_{2m}\vee O_{2n}$ has size $4mn+2m-1$, we define a bijection $f: E(P_{2m}\vee O_{2n})\to [1,4mn+2m-1]$ such that $f(u_iu_{i+1}) = i/2$ for even $i$, $f(u_iu_{i+1}) = 2m - (i+1)/2$ for odd $i$, and

\begin{enumerate}[(1)]
  \item $f(u_{2m-1}v_1) = 2m$ and $f(u_{2i-1}v_1) = 4m - i$ for $1\le i\le m-1$;
  \item $f(u_{2m}v_1) = 4mn+2m-1$ and $f(u_{2i}v_1) = 4mn + m - 1 - i$ for $1\le i\le m-1$; 
  \item $f(u_{2m-1}v_2) = 3m$ and $f(u_{2i-1}v_2) = 3m-i$ for $1\le i\le m-1$;
  \item $f(u_{2i}v_2) = 4mn + m - 2 + i$ for $1\le i\le m$;
  \item $f(u_{2i-1}v_3) = 4mn - m + i - 1$ for $1\le i\le m$;
  \item $f(u_{2m}v_3) = 4m$, $f(u_{2i}v_3) = 4m+i$ for $1\le i\le m-1$;
  \item $f(u_{2i-1}v_{2j}) = (2j+1)m - 1 + i$ for $1\le i\le m, 2\le j\le n$;
  \item $f(u_{2i}v_{2j}) = (4n+3-2j)m - i$ for $1\le i\le m, 2\le j\le n$;
  \item $f(u_{2i-1}v_{2j-1}) = (4n+4-2j)m - i$ for $1\le i\le m, 3\le j\le n$;
  \item $f(u_{2i}v_{2j-1}) = 2jm - 1 + i$ for $1\le i\le m, 3\le j\le n$.
\end{enumerate}

\nt It is easy to verify that $f^+(v_j) = m(4mn + 4m - 1)$ for $1\le j\le n$. Note that the values contributed by $P_{2m}$ to $f^+(u_{2i-1})$ for $1\le i\le m$, to $f^+(u_{2i})$ for $1\le i\le m-1$ and to $f^+(u_{2m})$ are $2m-1$, $2m$ and $m$ respectively. 
Thus, it is routine to check that 
\begin{eqnarray*}
f^+(u_{2m-1}) &=& (2m-1) + 2m + 3m + (4mn-1) + (6m-1) + \sum^n_{j=3}[m(4n+5)-1]\\
 &=& m(4n^2+n+3)-n-1.
\end{eqnarray*}

\nt Moreover, for $1\le i\le m-1$, 
\begin{eqnarray*}
f^+(u_{2i-1}) &=& (2m-1) + (4m-i) + (3m-i) + (4mn-m+i-1) + (5m-1+i) + \sum^n_{j=3}[m(4n+5)-1]\\
 &=& m(4n^2+n+3)-n-1.
\end{eqnarray*}

\nt Similarly, 
\begin{eqnarray*} 
f^+(u_{2m}) &=& m + (4mn+2m-1) + (4mn+2m-2) + 4m + (4n-2)m + \sum^n_{j=3}[m(4n+3)-1] \\
&=& m(4n^2 + 7n+1)-n-1. 
\end{eqnarray*} 

\nt Moreover, for $1\le i\le m-1$, 

\begin{eqnarray*}
f^+(u_{2i}) &=& 2m + (4mn+m-1-i) + (4mn+m-2+i) + (4m+i) + [(4n-1)m-i] + \sum^n_{j=3}[m(4n+3)-1] \\
&=& m(4n^2 + 7n+1)-n-1.
\end{eqnarray*} 

\ms\nt Clearly, $f^+(u_{2i-1}) < f^+(u_{2i})$ for $1\le i\le m$. Now, $f^+(v_j) - f^+(u_{2i-1}) = m(4mn+4m-4n^2-n-4) + n + 1 > 0$ if $m\ge n$. Otherwise, $m\le n-1$ and $f^+(v_j) - f^+(u_{2i-1}) \le (n-1)[4(n-1)n + 4(n-1)-4n^2-n-4]+n+1 = (n-1)(-n-8)+n+1 < 0$. Similarly, $f^+(v_j)\ne f^+(u_{2i})$. Thus, $f$ is a local antimagic labeling that induces $3$ distinct vertex colors. Thus, $\chi_{la}(P_{2m} \vee O_{2n}) \le 3$. Since $\chi_{la}(P_{2m}\vee O_{2n}) \ge \chi(P_{2m}\vee O_{2n}) = 3$, the theorem holds.
\end{proof} 

\begin{example} The labeling matrix of $P_6\vee O_8$ under $f$ is given below. The edge labels of $P_6$ are $5,1,4,2,3$ consecutively.
\[\begin{array}{c|*{8}{c}|c|c}
& v_1 & v_2 & v_3 & v_4 & v_5 & v_6 & v_7 & v_8 & from  \,\,P_6\,\, edges & f^+(u_i)\\\hline
u_1 & 11 & 8 & 45 & 15 & 41 & 21 & 35 & 27 & 5 & 208\\
u_3 & 10 & 7 & 46 & 16 & 40 & 22 & 34 & 28 & 5 & 208\\
u_5 & 6 & 9 & 47 & 17 & 39 & 23 & 33 & 29 & 5 & 208\\\hline
u_2 & 49 & 50 & 13 & 44 & 18 & 38 & 24 & 32 & 6 & 274 \\
u_4 & 48 & 51 & 14 & 43 & 19 & 37 & 25 & 31 & 6 & 274 \\
u_6 & 53 & 52 & 12 & 42 & 20 & 36 & 26 & 30 & 3 & 274 \\\hline
f^+(v_j) & 177 & 177 & 177 & 177 & 177 & 177 & 177  & 177 &
\end{array}\]
\end{example}

\ms
\begin{theorem}\label{thm-P2mVO2n-1} For $m, n\ge 1$, $\chi_{la}(P_{2m}\vee O_{2n-1}) = 3$ except that $\chi_{la}(P_4\vee O_1) = 4$.   \end{theorem}  

\begin{proof} When $m=1$, $P_2\vee O_{2n-1} = K_{1,1,2n-1}$. In~\cite[Theorem 2.7]{LNS-MJM}, the authors proved that $\chi_{la}(K_{1,1,2n-1}) = 3$.  When $n=1$, $P_{2m}\vee O_1 = F_{2m}$ the fan graph of order $2m+1$. In~\cite[Theorems 3.5 \& 3.6]{LSN}, the authors proved that $\chi_{la}(P_4\vee O_1) = 4$ and $\chi_{la}(P_{2m}\vee O_1) = 3$ for $m\ge 3$. 

\ms\nt We now consider $m, n\ge 2$. Since $P_{2m}\vee O_{2n-1}$ has size $4mn-1$, we define a bijection $f: E(P_{2m}\vee O_{2n-1})\to [1,4mn-1]$ such that $f(u_iu_{i+1}) = i/2$ for even $i$, $f(u_iu_{i+1}) = 2m - (i+1)/2$ for odd $i$, and

\begin{enumerate}[(1)]
  \item $f(u_{2m-1}v_1) = 2m$ and $f(u_{2i-1}v_1) = 4m - i$ for $1\le i\le m-1$;
  \item $f(u_{2m}v_1) = 4mn-1$ and $f(u_{2i}v_1) = 4mn - 2m + i - 1$ for $1\le i\le m-1$; 
  \item $f(u_{2m-1}v_2) = 3m$ and $f(u_{2i-1}v_2) = 3m-i$ for $1\le i\le m-1$;
  \item $f(u_{2i}v_2) = 4mn - m - 2 + i$ for $1\le i\le m$;
  \item $f(u_{2i-1}v_{2n-1}) = 2mn + 2(i-1)$ for $1\le i\le m$;
  \item $f(u_{2i}v_{2n-1}) = 2mn + 2m + 1 - 2i$ for $1\le i\le m$;
  \item $f(u_{2i-1}v_{2j-1}) = 4mn + 2m - 2jm - i$ for $1\le i\le m$ and $2\le j\le n-1$;
  \item $f(u_{2i}v_{2j-1}) = 2jm + i-1$ for $1\le i\le m$ and $2\le j\le n-1$;
  \item $f(u_{2i-1}v_{2j}) = (2j+1)m + i - 1$ for $1\le i\le m$ and $2\le j\le n-1$;
  \item $f(u_{2i}v_{2j}) = 4mn + m - 2jm - i$ for $1\le i\le m$ and $2\le j\le n-1$.
\end{enumerate}

\nt It is easy to verify that $f^+(v_j) = m(4mn + 2m - 1)$ for $1\le i\le 2n-1$. Note that the values contributed by $P_{2m}$ to $f^+(u_{2i-1})$ for $1\le i\le m$, to $f^+(u_{2i})$ for $1\le i\le m-1$ and to $f^+(u_{2m})$ are $2m-1$, $2m$ and $m$ respectively. Suppose $n=2$, we have $f^+(v_j) = m(10m-1)$ for $j=1,2,3$. It is routine to check that $f^+(u_{2i-1}) = 13m-3$ for $1\le i\le m$ and $f^+(u_{2i})= 21m-2$ for $1\le i\le m$. Clearly, $f$ is a local antimagic labeling that induces $3$ distinct vertex labels. Thus, $\chi_{la}(P_{2m} \vee O_3)\le 3$. 

\ms\nt For $n\ge 3$, it is routine to check that

\begin{eqnarray*}
f^+(u_{2m-1}) &=& (2m-1) + 2m + 3m + (2mn+2m-2) + \sum^{n-1}_{j=2}[m(4n+3)-1]\\
 &=& m(4n^2-3n+3)-n-1.
\end{eqnarray*}

\nt Moreover, for $1\le i\le m-1$, 
\begin{eqnarray*}
f^+(u_{2i-1}) &=& (2m-1) + (4m-i) + (3m-i) + (2mn + 2i-2) + \sum^{n-1}_{j=2}[m(4n+3)-1]\\
 &=& m(4n^2-3n+3)-n-1.
\end{eqnarray*}

\nt Similarly, 
\begin{eqnarray*} 
f^+(u_{2m}) &=& m + (4mn-1) + (4mn-2) + (2mn+1) + \sum^{n-1}_{j=2}[m(4n+1)-1] \\
&=& m(4n^2 + 3n - 1)-n. 
\end{eqnarray*} 

\nt Moreover, for $1\le i\le m-1$, 

\begin{eqnarray*}
f^+(u_{2i}) &=& 2m + (4mn-2m+i-1) + (4mn-m-2+i) + (2mn+2m+1-2i) + \sum^{n-1}_{j=2}[m(4n+1)-1] \\
&=& m(4n^2 + 3n - 1)-n.
\end{eqnarray*} 

\nt Clearly, $f^+(u_{2i-1}) < f^+(u_{2i})$ for $1\le i\le m$. Now, $f^+(v_j) - f^+(u_{2i-1}) = m(4mn + 2m - 4n^2 + 3n - 4) + n + 1 > 0$ if $m\ge n$. If $m = n-1$, then $f^+(v_j) - f^+(u_{2i-1}) = (n-1)[4(n-1)n + 2(n-1) - 4n^2 + 3n - 4] + n + 1 = n^2-6n+7\ne 0$. Otherwise, $m\le n-2$ and $f^+(v_j) - f^+(u_{2i-1})\le (n-2)[4(n-2)n + 2(n-2) - 4n^2 + 3n - 4] + n + 1 = (n-2)(-3n-8)+n+1 < 0$. Similarly, $f^+(v_j) - f^+(u_{2i})\ne 0$. Thus, $f$ is a local antimagic labeling that induces $3$ distinct vertex colors. Thus, $\chi_{la}(P_{2m} \vee O_{2n-1}) \le 3$. Since $\chi_{la}(P_{2m}\vee O_{2n-1}) \ge \chi(P_{2m}\vee O_{2n-1}) = 3$, the theorem holds.
\end{proof} 

\begin{example}\label{Eg-P6VO5}  The labeling matrix of $P_6\vee O_5$ under $f$ is given below. The edge labels of $P_6$ are $5,1,4,2,3$ consecutively.
\[\begin{array}{c|*{5}{c}|c|c}
& v_1 & v_2 & v_3 & v_4 & v_5 & from\,\, P_6 \,\,edges & f^+(u_i)\\\hline
u_1 & 11 & 8 & 29 & 15 & 18 & 5 & 86\\
u_3 & 10 & 7 & 28 & 16 & 20 & 5 & 86\\
u_5 & 6 & 9 & 27 & 17 & 22 & 5 & 86\\\hline
u_2 & 30 & 32 & 12 & 26 & 23 & 6 & 129 \\
u_4 & 31 & 33 & 13 & 25 & 21 & 6 & 129 \\
u_6 & 35 & 34 & 14 & 24 & 19 & 3 & 129 \\\hline
f^+(v_j) & 123 & 123 & 123 & 123  & 123 &
\end{array}\]
\end{example}

\nt We now give an ad hoc example to show that $\chi_{la}(P_7 \vee O_3)=3$ as in the labeling matrix below.

\begin{example} The edge labels of $P_7$ are $4,1,5,2,6,3$ consecutively.
\[\begin{array}{c|*{3}{c}|c|c}
& v_1 & v_2 & v_3 & from \,\,P_7\,\, edges & f^+(u_i)\\\hline
u_1 & 12 & 14 & 21 & 4 & 51\\
u_3 & 15 & 10 & 20 & 6 & 51\\
u_5 & 19 & 16 & 8 & 8 & 51 \\
u_7 & 13 & 17 & 18 & 8 & 51 \\\hline
u_2 & 27 & 11 & 22 & 5 & 65 \\
u_4 & 9 & 26 & 23 & 7 & 65 \\
u_6 & 24 & 25 & 7 & 9 & 65 \\\hline
f^+(v_j) & 119 & 119 & 119 & & 
\end{array}\]
\end{example}

\nt We now consider more join graphs.

\begin{theorem}\label{thm-P2mVC2n-1} For $m\ge 1$, $n\ge 2$, $\chi_{la}(P_{2m}\vee C_{2n-1}) = 5$. \end{theorem} 

\begin{proof} We note that $\chi(P_{2m}\vee C_{2n-1}) = 5$. Let $C_{2n-1} = (v_1,v_2,v_3\cdots, v_{2n-1},v_1)$. When $m=1$, $P_2\vee C_{2n-1}$ has size $6n-2$. We define a bijection $g: E(P_2\vee C_{2n-1})\to [1,6n-2]$ such that $g(u_1v_j) = j$, $g(u_2v_j) = 4n-1-j$ for $1\le j\le 2n-1$ and $g(u_1u_2) = 4n-1$. Moreover, $g(v_jv_{j+1}) = 4n-1 + j/2$ for even $j$, and $g(v_jv_{j+1}) = 6n-2-(j-1)/2$. 

\ms\nt Clearly, $g^+(u_1) = 2n^2-n$, $g^+(u_2) = 6n^2-5n+1$. Moreover, $g^+(v_1) = 15n-4$, $g^+(v_{2j-1}) = 14n-4$ for $2\le j\le n$ and $g^+(v_{2j}) = 14n-3$ for $1\le j\le n$. Since the induced vertex colors are distinct, $g$ is a local antimagic labeling so that $\chi_{la}(P_2\vee C_{2n-1})\le 5$.

\ms\nt Consider $m\ge 2$. Let $f$ be the local antimagic labeling defined in the proof of Theorem~\ref{thm-P2mVO2n-1}. Since $P_{2m}\vee C_{2n-1}$ has size $4mn + 2n - 2$, we define a bijection $g: E(P_{2m}\vee C_{2n-1}) \to [1, 4mn+2n-2]$ such that $g(e) = f(e)$ if $e\in E(P_{2m}\vee O_{2n-1})$, and that $g(v_jv_{j+1}) = 4mn-1 + j/2$ for even $j$, and $g(v_jv_{j+1}) = 4mn+2n-2-(j-1)/2$ for odd $j$.

\ms\nt Clearly, $g^+(u_{2i}) = f^+(u_{2i})$, and $g^+(u_{2i-1}) = f^+(u_{2i-1})$ for $1\le i\le m$. Moreover, $g^+(v_1) = f^+(v_1) + 8mn + 3n - 3$, $g^+(v_{2j-1}) = f^+(v_{2j-1}) + 8mn + 2n - 3$ for $2\le j\le n$, and $g^+(v_{2j}) = f^+(v_{2j}) + 8mn + 2n-2$ for $1\le j\le n$.

\ms\nt When $n=2$, the $5$ induced vertex colors are $g^+(u_{2i-1})=13m-3$ and $g^+(u_{2i})=21m-2$ for $1\le i\le m$, and $g^+(v_1)=10m^2+15m+3$, $g^+(v_2)=10m^2+15m+2$, $g^+(v_3)=10m^2+15m+3$ which are all distinct. 

\ms\nt When $n\ge 3$, the $5$ induced vertex colors are 
\begin{enumerate}[(1)]
  \item $g^+(u_{2i-1})=m(4n^2-3n+3)-n-1$ for $1\le i\le m$;
  \item $g^+(u_{2i})=m(4n^2+3n-1)-n$ for $1\le i\le m$;
  \item $g^+(v_1)=m(4mn+2m+8n-1)+3n-3$,
  \item $g^+(v_{2j-1})=m(4mn+2m+8n-1)+2n-3$ for $2\le j\le n$;
  \item $g^+(v_{2j})=m(4mn+2m+8n-1)+2n-2$ for $1\le j\le n-1$. 
\end{enumerate} 
\nt As in the proof of Theorem~\ref{thm-P2mVO2n-1}, it is easy to verify that $(1) \ne (2)\ne (3) \ne (4) \ne (4)$.  Consequently, $g$ is a local antimagic labeling so that $\chi_{la}(P_{2m}\vee C_{2n-1})\le 5$. Since $\chi_{la}(P_{2m}\vee C_{2n-1})\ge \chi(P_{2m}\vee C_{2n-1})=5$, the theorem holds.
\end{proof} 

\begin{example} To get a labeling of $P_6\vee C_5$, we add edges to the vertices $v_j, 1\le j\le 5$, of $P_6\vee O_5$ to form $C_5=v_1v_2v_3v_4v_5v_1$ and label the edges by $40,36,39,37,38$ consecutively. By referring to Example~\ref{Eg-P6VO5}, it is easy to verify that the induced vertex labels of $u_i, 1\le i\le 6$ remain unchanged while the induced vertex labels of $v_1$ to $v_5$ are $201,199,198,199,198$ respectively.
\end{example}

\begin{theorem}\label{thm-P2mVK2n} For $m\ge 1$, $n\ge 1$, $\chi_{la}(P_{2m} \vee K_{2n})= 2n+2$. \end{theorem}
  
\begin{proof} Let $G=P_{2m} \vee K_{2n}$ be obtained from $P_{2m}\vee O_{2n}$ by adding $n(2n-1)$ edges joining $v_i,v_j$ for $1\le i < j \le 2n$. If $m=1$, then $G = K_{2n+2}$ with $\chi_{la}(G)=2n+2$. We now consider $m\ge 2$. 
Suppose $n=1$, then $|E(G)|=6m$. By making minor changes to the labeling function $f$ in~\cite[Theorem 2.4]{YBYL}, we define an edge labeling $g:E(P_{2m}\vee K_{2})\to [1,6m]$ such that 
$$g(u_iu_{i+1}) = \begin{cases} 2m - (i+1)/2 & \mbox{ if $i$ is odd,} \\ i/2 &\mbox{ if $i$ is even.}\end{cases}$$
$$g(u_iv_1) = \begin{cases} 2m + (i-1)/2 & \mbox{ if $i$ is odd,} \\ 6m - (i+2)/2 & \mbox{ if $i$ is even, $i\ne  2m$}. \end{cases}$$
$$g(u_iv_2) = \begin{cases} 5m - (i+1)/2 & \mbox{ if $i$ id odd,} \\ 3m + (i-2)/2 & \mbox{ if $i$ is even, $i\ne  2m$.} \end{cases}$$
$$g(u_{2m}v_1) = 4m-1, g(u_{2m}v_2) = 6m-1, g(v_1v_2) = 6m$$

\nt Thus, the induced vertex labels are $g^+(u_i) = f^+(u_i) = 9m-2$ for odd $i$, $g^+(u_i) = f^+(u_i) = 11m-2$ for even $i$, $g^+(v_1) = f^+(v_1)-(6m-1)+(4m-1)+6m = 8m^2+3m$, $g^+(v_2) = f^+(v_2) - (4m-1) + (6m-1) + 6m = 8m^2+7m$. Clearly, all the four induced vertex colors are distinct. 

\ms\nt Suppose $n\ge 2$. Let $h: E(K_{2n})\to [1,n(2n-1)]$ be a local antimagic labeling of $K_{2n}.$ Note that $h^+(v_j)$ are distinct for $1\leq j\leq 2n.$  Let $f$ be the local antimagic labeling of $P_{2m}\vee O_{2n}$ as defined in the proof of Theorem~\ref{thm-P2mVO2n}. Define $g: E(P_{2m}\vee K_{2n})\to [1,4mn+2m-1+n(2n-1)] $ such that $g(e)=f(e)$ for $e\in E(P_{2m}\vee O_{2n})$ and $g(e)=h(e)+4mn+2m-1$ for $e\in E(K_{2n}).$

\ms\nt Thus, for $1\le i\le m$,  $g^+(u_{2i-1}) = f^+(u_{2i-1}) = m(4n^2+n+3)-n-1$ and $g^+(u_{2i}) = f^+(u_{2i}) = m(4n^2 + 7n+1)-n-1$. Moreover, for $1\le j\le 2n$, $g^+(v_{j}) = h^+(v_j) + f^+(v_j) + (2n-1)(4mn+2m-1) = h^+(v_j) + m(4mn+4m-1) + (2n-1)(4mn+2m-1)$. Since $h^+(v_j)$ are all distinct, we have $g^+(v_j)$ are all distinct too.

\ms\nt Clearly, $g^+(u_{2i-1}) < g^+(u_{2i})$ for $1\le i\le m$. Now, $g^+(v_j) - g^+(u_{2i-1}) = h^+(v_j) + m(4mn+4m-1) + (2n-1)(4mn+2m-1)-[m(4n^2+n+3)-n-1]= h^+(v_j) + m(4mn+4n^2+4m-6-n)-n+2 > 0$ for all $m,n\ge 2$.  Thus $g^+(v_j)\ne g^+(u_{2i-1})$. Similarly, $g^+(v_j)\ne g^+(u_{2i})$. Thus, $g$ is a local antimagic labeling that induces $2n+2$ distinct vertex colors. Thus, $\chi_{la}(P_{2m} \vee K_{2n}) \le 2n+2$. Since $\chi_{la}(P_{2m}\vee K_{2n}) \ge \chi(P_{2m}\vee K_{2n}) = 2n+2$, the theorem holds.
\end{proof}

\begin{theorem}\label{thm-P2mVK2n-1} For $m\ge 1$, $n\ge 1$, $\chi_{la}(P_{2m} \vee K_{2n-1})= 2n+1$ except that $\chi_{la}(P_4\vee K_1) = 4$. \end{theorem}

\begin{proof} Let $G=P_{2m} \vee K_{2n-1}$ be obtained from $P_{2m}\vee O_{2n-1}$ by adding $(n-1)(2n-1)$ edges joining $v_i,v_j$ for $1\le i < j \le 2n-1$. If $m=1$, then $G = K_{2n+1}$ with $\chi_{la}(G)=2n+1$. We now consider $m\ge 2$. 
Suppose $n=1$, then $G=F_{2m}$, the fan graph of order $2m+1.$ In~\cite[Theorems 3.5 \& 3.6]{LSN}, the authors proved that $\chi_{la}(P_4\vee O_1) = 4$ and $\chi_{la}(P_{2m}\vee O_1) = 3$ for $m\ge 3$. When $n=2$, $K_3 = C_3$ and the result follows from Theorem~\ref{thm-P2mVC2n-1}.

\ms\nt We now consider $n\ge 3$. Let $h: E(K_{2n-1})\to [1,(n-1)(2n-1)]$ be a local antimagic labeling of $K_{2n-1}.$ Note that $h^+(v_j)$ are distinct for $1\leq j\leq 2n-1.$  Let $f$ be the local antimagic labeling of $P_{2m}\vee O_{2n-1}$ as defined in the proof of Theorem~\ref{thm-P2mVO2n-1}. Define $g: E(P_{2m}\vee K_{2n-1})\to [1,4mn-1+(n-1)(2n-1)] $ such that $g(e)=f(e)$ for $e\in E(P_{2m}\vee O_{2n-1})$ and $g(e)=h(e)+4mn-1$ for $e\in E(K_{2n-1}).$

\ms\nt Thus, for $1\le i\le m$,  $g^+(u_{2i-1}) = f^+(u_{2i-1}) = m(4n^2-3n+3)-n-1$ and $g^+(u_{2i}) = f^+(u_{2i}) = m(4n^2 + 3n-1)-n$. Moreover, for $1\le j\le 2n-1$, $g^+(v_{j}) = h^+(v_j) + f^+(v_j) + (2n-2)(4mn-1) = h^+(v_j) + m(4mn+2m-1)+(2n-2)(4mn-1)$. Since $h^+(v_j)$ are all distinct, we have $g^+(v_j)$ are all distinct too.

\ms\nt Clearly, $g^+(u_{2i-1}) < g^+(u_{2i})$ for $1\le i\le m$. Now, $g^+(v_j) - g^+(u_{2i-1}) = h^+(v_j) + m(4mn+2m-1)+(2n-2)(4mn-1)-m(4n^2-3n+3)-n-1 = h^+(v_j) + m(4mn+2m-4+4n^2-5n)-n+3$ for $m\ge 2, n\ge 3$. Thus $g^+(v_j)\ne g^+(u_{2i-1})$. Similarly, $g^+(v_j)\ne g^+(u_{2i})$. Thus, $g$ is a local antimagic labeling that induces $2n+1$ distinct vertex colors. Thus, $\chi_{la}(P_{2m} \vee K_{2n-1}) \le 2n+1$. Since $\chi_{la}(P_{2m}\vee K_{2n-1}) \ge \chi(P_{2m}\vee K_{2n-1}) = 2n+1$, the theorem holds.
\end{proof}

\nt In~\cite[Theorems 3.1 \& 3.2]{LSN}, the authors determined the exact value for $\chi_{la}(C_m \vee O_n)$ for $m\equiv n\pmod{2}$. 

\begin{theorem}\label{thm-C2mVO2n-1} For $m\ge 2,n\ge 1$, $\chi_{la}(C_{2m}\vee O_{2n-1}) = 3$.
\end{theorem} 

\begin{proof} Let $C_{2m} = (u_1,u_2,\cdots, u_{2m},u_1)$. Note that $C_{2m}\vee O_1 = W_{2m}$, the wheel graph of order $2m+1$. In~\cite[Theorems 2.14 \& 5]{Arumugam, LNS}, the authors proved that $\chi_{la}(W_{2m}) = 3$ for $m\ge 2$. We consider $n\ge 2$.

\ms\nt Since $C_{2m}\vee O_{2n-1}$ has size $4mn$, we define a bijection $g : E(C_{2m}\vee O_{2n-1}) \to [1,4mn]$ such that $g(u_1u_{2m}) = m+1$, $g(u_iu_{i+1}) = m - (i-1)/2$ for odd $1\le i\le 2m$, $g(u_iu_{i+1}) = m+1+i/2$ for even $2\le i\le 2m-2$, $g(u_iv_j) = f(u_iv_j) + 1$ for $1\le i\le m, 1\le j\le n$ where $f$ is the local antimagic labeling defined in the proof of Theorem~\ref{thm-P2mVO2n-1}. Note that the values contributed by $C_{2m}$ to $f^+(u_{2i-1})$ for $1\le i\le m$, to $f^+(u_{2i})$ for $1\le i\le m-1$ and to $f^+(u_{2m})$ are $2m+1$, $2m+2$ and $m+2$ respectively. There is an increment of $2$ compare to the original values contributed by $P_{2m}$. As such, $g^+(u_i) = f^+(u_i) + 2 + (2n-1) = f^+(u_i) + 2n+1$ for $1\le i\le 2m$ and $g^+(v_j) = f^+(v_j) + 2m$ for $1\le j\le 2n-1$.

\ms\nt When $n=2$, we now have $g^+(u_{2i-1}) = 13m-3 + 5 = 13m+2$ and $g^+(u_{2i}) = 21m-2+5=21m+3$ for $1\le i\le m$ whereas $g^+(v_j) = m(10m-1)+2m = m(10m+1)$. Clearly, $g$ is a local antimagic labeling that induces $3$ distinct vertex colors.

\ms\nt Consider $n\ge 3$. We now have $g^+(u_{2i-1}) = m(4n^2-3n+3)-n-1+(2n+1) = m(4n^2-3n+3) + n$ and $g^+(u_{2i}) = m(4n^2+3n-1)-n+(2n+1) = m(4n^2+3n-1)+n+1$ for $1\le i\le m$. Moreover, $g^+(v_j) = m(4mn+2m-1) + 2m = m(4mn + 2m + 1)$ for $1\le j\le 2n-1$. It is routine to check as in the proof of Theorem~\ref{thm-P2mVO2n-1} that the induced vertex colors are distinct. Thus, $g$ is a local antimagic labeling and $\chi_{la}(C_{2m}\vee O_{2n-1})\le 3$. Since $\chi_{la}(C_{2m}\vee O_{2n-1})\ge \chi(C_{2m}\vee O_{2n-1}) = 3$, the theorem holds.
\end{proof} 

\begin{example}\label{Eg-C6VO5}  The labeling matrix of $C_6\vee O_5$ under $g$ is given below. The edge labels of $C_6=u_1u_2u_3u_4u_5u_6u_1$ are $3,5,2,6,1,4$ consecutively.
\[\begin{array}{c|*{5}{c}|c|c}
& v_1 & v_2 & v_3 & v_4 & v_5 & \rm{from}  \,\,C_6\,\, \rm{ edges} & f^+(u_i)\\\hline
u_1 & 12 & 9 & 30 & 16 & 19 & 7 & 93\\
u_3 & 11 & 8 & 29 & 17 & 21 & 7 & 93\\
u_5 & 7 & 10 & 28 & 18 & 23 & 7 & 93\\\hline
u_2 & 31 & 33 & 13 & 27 & 24 & 8 & 136 \\
u_4 & 32 & 34 & 14 & 26 & 22 & 8 & 136 \\
u_6 & 36 & 35 & 15 & 25 & 20 & 5 & 136 \\\hline
f^+(v_j) & 129 & 129 & 129 & 129  & 129 &
\end{array}\]
\end{example}

\begin{theorem} For $n\ge 1$, $\chi_{la}(C_{2n+1}\vee O_{2n}) = 4$.  \end{theorem}

\begin{proof} Let $G = C_{2n+1}\vee O_{2n}$ with $V(G) = \{u_i, v_j\} \,|\, 1\le i\le 2n+1, 1\le j\le 2n\}$. Let $M$ be the $(2n+1)\times(2n+1)$ magic square obtained using the Siamese method~\cite{Krai}. Note that the $(i,n)$-th entry is $1+2(n+1)(i-1)$ for $i=1,2,\ldots,2n+1$. The row and column sum of $M$ is $K = (2n+1)[(2n+1)^2+1]/2$. Let $f : E(G) \to [1,(2n+1)^2]$ be an edge labeling of $G$. We now describe the procedure to label the edges of $G$ in the following 4 steps.

\begin{enumerate}[(1)]
  \item Delete the $n$-th column of $M$ to get a $(2n+1)\times2n$ matrix, $N$ that has $i$-th row sum $= K - 1-2(n+1)(i-1)$ for $i=1,\ldots,2n+1$.
  \item Move row $(2n+1)$ to row 1 and row $i$ to row $(i+1)$ for $i=1,3,\ldots,2n-1$ to get a new matrix $P$.
  \item Let $f(u_iv_j)$ be the $(i,j)$-entry of $P$. 
  \item Label the edges of $C_{2n+1} = u_1u_2u_3\cdots u_{2n-1}u_{2n}u_{2n+1}u_1$ by $f(u_{2i-1}u_{2i}) = 1+2(n+1)(i-1)$ for $i=1,2,3,\ldots,n+1$, where $u_{2n+2} = u_1$, and that $f(u_{2i}u_{2i+1}) = 1+2(n+1)(n+i)$ for $i=1,2,3,\ldots,n$. 
\end{enumerate}

\nt We observe that 
\begin{enumerate}[(a)]
  \item $f$ is bijective with $f^+(v_j) = K$ for $j=1,\ldots,2n$.
  \item The edge labels of $C_{2n+1}$ contribute $2+2n(n+1)+(i-1)(2n+2)$ to $u_i$ for $i=1,2,\ldots,2n+1$.
  \item Matrix $P$ has row 1 sum $= K - 1 - 4n(n+1)$. For $i=1,2,\ldots,n$, row $(2i+1)$ has sum $= K - 1 - 2(n+1)(2i-2)$ whereas row $(2i)$ has sum $= K - 1 - 2(n+1)(2i-1)$.
\end{enumerate} 

\nt We now have 
\begin{enumerate}[(i)]
  \item $f^+(v_j) = K$ for $j=1,\ldots, 2n$, and $f^+(u_1) = K-1-4n(n+1) + 2 + 2n(n+1) = K + 1 - 2n(n+1)$,
  \item $f^+(u_{2i+1}) = K-1-2(n+1)(2i-2) + 2+2n(n+1)+(2i)(2n+2) = K+1+(4+2n)(n+1)$ for $i=1,\ldots,n$,
  \item $f^+(u_{2i}) = K-1-2(n+1)(2i-1) + 2+2n(n+1)+(2i-1)(2n+2) = K+1+2n(n+1)$ for $i=1,\ldots,n$.
\end{enumerate} 

\nt Clearly, $f$ is a local antimagic $4$-labeling of $G$. Thus, $\chi_{la}(G)\le 4$. Since $\chi_{la}(G) \ge \chi(G) = 4$, the theorem holds.
\end{proof}

\begin{example}\label{Eg-C7VO6}  The labeling matrix of $C_7\vee O_6$ under $f$ is given below. The edge labels of $C_7=u_1u_2u_3u_4u_5u_6u_7u_1$ are $1,33,9,41,17,49,25$ consecutively.
\[\begin{array}{c|*{6}{c}|c|c}
& v_1 & v_2 & v_3 & v_4 & v_5 & v_6 & \rm{from}  \,\,C_7\,\, \rm{ edges} & f^+(u_i)\\\hline
u_1 & 22 & 31 & 40 & 2 & 11 & 20 & 26 & 152\\
u_2 & 38 & 47 & 7 & 18 & 27 & 29 & 34 & 200\\
u_3 & 30 & 39 & 48 & 10 & 19 & 28 & 42 & 216\\
u_4 & 5 & 14 & 16 & 34 & 36 & 45 & 50 & 200 \\
u_5 & 46 & 6 & 8 & 26 & 35 & 37 & 58 & 216 \\
u_6 & 21 & 23 & 32 & 43 & 3 & 12 & 66 & 200 \\
u_7 & 13 & 15 & 24 & 42 & 44 & 4 & 74 & 216 \\\hline
f^+(v_j) & 175 & 175 & 175 & 175  & 175 & 175 
\end{array}\]
\end{example}

\nt It is easy to verify the conditions of Lemmas~\ref{lem-nonreg} and~\ref{lem-G-e}. Thus, we have the following corollary.

\begin{corollary} Suppose $e\ne E(C_{2n+1})$. For $n\ge 1$, $\chi_{la}((C_{2n+1}\vee O_{2n}) - e) = 4$. \end{corollary}

\nt In~\cite[Theorems 3.5 \& 3.6]{LSN}, the authors determined $\chi_{la}((C_{2m}\vee O_1)-e)$. We now have the following general result.

\begin{theorem}\label{thm-C2mVO2n-1-e} For $m,n\ge 2$, $\chi_{la}((C_{2m}\vee O_{2n-1})-e) = 3$. \end{theorem} 

\begin{proof} Let $G = C_{2m} \vee O_{2n-1}$ defined in the proof of Theorem~\ref{thm-C2mVO2n-1} and $H=G - e$. We consider the two cases.

\ms\nt {\bf Case (a).} $e \in E(C_{2m})$. Let  $g$ be the antimagic labeling defined in the proof of Theorem~\ref{thm-C2mVO2n-1}. Without loss of generality, we may assume $e=u_{2m-1}u_{2m}$. Now, $g(e)=1$. Moreover, $g$ induces a 3-independent partitions such that all vertices in the same independent set have the same degree. Precisely, $deg(v_j) = 2m$ for $1\le j\le 2n-1$, $deg(u_{2i-1}) = 2n+1$ for $1\le i\le m$ and $ deg(u_{2i})= 2n+1$ for $1\le i\le m$. Thus, 

\begin{enumerate}[(1)]
  \item $g^+(u_{2i-1}) - deg(u_{2i-1}) = 4mn^2-3nm+3m-n-1 $ for $1\le i\le m$,
  \item $g^+(u_{2i}) - deg(u_{2i}) = 4mn^2+3nm-m-n $ for $1\le i\le m$,
  \item $g^+(v_j) - deg(v_j) = 4m^2n+2m^2-m$ for $1\le j\le 2n-1$.
\end{enumerate} 

\nt Clearly, $(1)-(2)= -6mn+4m-1 \ne 0$. Thus $(1)\ne(2)$. Now $(1)-(3)= 4mn(n-m)-3nm-2m^2+4m-n-1 > 0$, if $m \le n-2$ and $(1)-(3)\le -n^2+6n-7 < 0$ if $m \ge n-1$. Thus $(1) \ne (3)$. Similarly, $(2)-(3)= 4mn(n-m)+3nm-2m^2-n > 0$, if $m \le n$ and $(2)-(3) \le -3n^2-6n-2 < 0$ if $m \ge n+1$. Thus $(1) \ne (2) \ne (3)$.

\ms\nt We observe that the conditions $(i)$ and $(ii)$ of Lemma~\ref{lem-G-e} are satisfied and so $\chi(H) \le \chi_{la}(H) \le 3.$ 

\ms\nt {\bf Case (b).} $e \notin E(C_{2m})$. Let $g$ be the local antimagic labeling of $G$ defined in the proof of Theorem~\ref{thm-C2mVO2n-1} such that $g(u_{2m}v_1) = 4mn$. Clearly for any $u,v \in V(G),$ $g^+(u)=g^+(v)$ implies that $deg(u)=deg(v).$ Also if $g^+(u)\neq g^+(v)$ then $(4mn+1)(deg(u)-deg(v))\neq g^+(u) - g^+(v).$  Therefore by Lemma~\ref{lem-nonreg}, $h=4mn+1-g$ is also a local antimagic labeling of $G$ such that $h(u_{2m}v_1) = 1$ and that
\begin{eqnarray*}
h^+(u_{2i-1}) &=& deg(u_{2i-1}) (4mn+1) - g^+(u_{2i-1})\\
 &=& 4mn^2 + 7mn + n -3m + 1.
\end{eqnarray*}
\begin{eqnarray*} 
h^+(u_{2i}) &=& deg(u_{2i}) (4mn+1) - g^+(u_{2i}) \\
&=& 4mn^2 + mn + n + m. 
\end{eqnarray*} 
\begin{eqnarray*}
h^+(v_j) &=&  deg(v_j) (4mn+1) - g^+(v_j)  \\
&=& 4m^2n - 2m^2 + m.
\end{eqnarray*} 

\nt Hence $h$ induces a 3-independent partitions such that all vertices in the same independent set have the same degree. Precisely, $deg(v_j) = 2m$ for $1\le j\le 2n-1$, $deg(u_{2i-1}) = 2n+1$ for $1\le i\le m$ and $ deg(u_{2i})= 2n+1$ for $1\le i\le m$. Without loss of generality, let $e=u_{2m}v_1$. Now, 
\begin{enumerate}[(1)]
  \item $h^+(u_{2i-1}) - deg(u_{2i-1}) = 4mn^2+7nm-3m-n $ for $1\le i\le m$,
  \item $h^+(u_{2i}) - deg(u_{2i}) = 4mn^2+nm+m-n-1 $ for $1\le i\le m$,
  \item $h^+(v_j) - deg(v_j) = 4m^2n-2m^2-m$ for $1\le j\le 2n-1$.
\end{enumerate} 

\nt Consider $(1)-(2)= 6mn-4m+1 \ne 0$, thus $(1)\ne(2)$. Now $(1)-(3)= 4mn(n-m)+7nm+2m^2-2m-n > 0$, if $m \le n+2$ and $(1)-(3)= -3n^2-6n+12 < 0$ if $m \ge n+3$. Thus $(1) \ne (3)$. Similarly, $(2)\ne (3)$. By Lemma~\ref{lem-G-e}, we have $\chi(H) \le \chi_{la}(H) \le 3.$ Since $\chi(H) = 3$, the theorem holds. \end{proof}

\nt In~\cite[Theorem 3.8]{LSN}, the authors proved that $\chi_{la}(C_{2m-1}\vee C_{2n-1}) = 6$ for $m,n\ge 2$.

\begin{theorem}\label{thm-C2mVC2n-1} For $m,n\ge 2$, $\chi_{la}(C_{2m}\vee C_{2n-1}) = 5$.
\end{theorem} 

\begin{proof} Let $C_{2n-1} = (v_1,v_2,v_3,\cdots, v_{2n-1},v_1)$ and $g$ be the local antimagic labeling of $C_{2m}\vee O_{2n-1}$ as defined in the proof of Theorem~\ref{thm-C2mVO2n-1}. Since $C_{2m}\vee C_{2n-1}$ has size $4mn + 2n - 1$, we define a bijection $h : E(C_{2m}\vee C_{2n-1}) \to [1, 4mn+2n-1]$ such that $h(e) = g(e)$ for $e\in E(C_{2m}\vee O_{2n-1})$ whereas $h(v_jv_{j+1}) = 4mn + j/2$ for even $j$ and $h(v_jv_{j+1}) = 4mn+2n-1 - (j-1)/2$ for odd $j$.  Note that the values contributed by $C_{2n-1}$ to $h^+(v_1)$, $h^+(v_{2j-1})$ for $2\le j\le n$ and $h^+(v_{2j})$ for $1\le j\le n$ are $8mn+3n-1$, $8mn+2n-1$ and $8mn+2n$ respectively. Obviously, $h^+(u_i) = g^+(u_i)$ for $1\le i\le 2m$, $h^+(v_1) = g^+(v_1) + 8mn + (3n-1)$, $h^+(v_{2j-1}) = g^+(v_{2j-1}) + 8mn + (2n-1)$ for $2\le j\le n$, and $h^+(v_{2j}) = g^+(v_{2j}) + 8mn + 2n$ for $1\le j\le n$. 

\ms\nt When $n=2$, we have $h^+(u_{2i-1}) = 13m+2$ and $h^+(u_{2i}) = 21m+3$ for $1\le i\le m$ whereas $h^+(v_1) = m(10m+1) + 8m+5 = 10m^2 + 9m + 5$, $h^+(v_2) = m(10m+1) + 8m+4 = 10m^2+9m+4$ and $h^+(v_3) =  m(10m+1) + 8m+3 = 10m^2+9m+3$.

\ms\nt Consider $n\ge 3$. We have $h^+(u_{2i-1}) = m(4n^2-3n+3) + n$ and $h^+(u_{2i}) = m(4n^2+3n-1)+n+1$ for $1\le i\le m$. Moreover, $h^+(v_1) =  m(4mn + 2m + 1) + 8mn + (3n-1) = m(4mn + 2m + 8n +1)+3n-1$, $h^+(v_{2j-1}) = m(4mn + 2m + 1) + 8mn + (2n-1) = m(4mn + 2m + 8n +1) + 2n-1$ for $2\le j\le n$, and $h^+(v_{2j}) = m(4mn + 2m + 1) + 8mn + 2n = m(4mn + 2m + 8n +1) + 2n$ for $1\le j\le n$. It is routine to check as in the proof of Theorem~\ref{thm-P2mVO2n-1} that the induced vertex colors are distinct. Thus, $g$ is a local antimagic labeling and $\chi_{la}(C_{2m}\vee C_{2n-1})\le 5$. Since $\chi_{la}(C_{2m}\vee C_{2n-1})\ge \chi(C_{2m}\vee C_{2n-1}) = 5$, the theorem holds. 
\end{proof} 

\begin{example}\label{Eg-C6VC5} To get a labeling of $C_6\vee C_5$, we add edges to the vertices $v_j, 1\le j\le 5$, of $C_6\vee O_5$ to form $C_5=(v_1,v_2,v_3,v_4,v_5,v_1)$ and label the edges by $41,37,40,38,39$ consecutively. By referring to Example~\ref{Eg-C6VO5}, it is easy to verify that the induced vertex labels of $u_i, 1\le i\le 6$ remain unchanged while the induced vertex labels of $v_1$ to $v_5$ are $209,207,206,207,206$ respectively.
\end{example}

\begin{theorem}\label{thm-C2mVC2n-1-e} For $m,n\ge 2$, $\chi_{la}((C_{2m}\vee C_{2n-1})-e) = 5$ for $e\in E(C_{2m})$.
 \end{theorem}

\begin{proof}  Suppose $G = C_{2m} \vee C_{2n-1}$ and $H=G-e$ for $e\in E(C_{2m})$. Let $h$ be the local antimagic labeling defined in the proof of Theorem~\ref{thm-C2mVC2n-1} that induces a 5-independent partitions such that all vertices in the same independent set have the same degree. Precisely, $deg(v_j) = 2m+2$ for $1\le j\le 2n-1$, $deg(u_{2i-1}) = deg(u_{2i})= 2n+1$ for $1\le i\le m$. Thus, 

\begin{enumerate}[(1)]
  \item $h^+(u_{2i-1}) - deg(u_{2i-1}) = m(4n^2-3n+3)-n-1$ for $1\le i\le m$, 
  \item $h^+(u_{2i}) - deg(u_{2i}) = m(4n^2 + 3n-1)-n $ for $1\le i\le m-1$,
  \item $h^+(v_1) - deg(u_1) = m(4mn+2m+8n+1)+3n-2m-3 $, 
  \item $h^+(v_{2j-1}) - deg(v_j) =  m(4mn+2m+8n+1)+2n-2m-3 $ for $2\le i\le n$,  
  \item $h^+(v_{2j}) - deg(u_{2m}) = m(4mn+2m+8n+1)+2n-2m-2 $ for $1\le i\le n$. 
\end{enumerate} 

\nt Clearly $(3) >(5) > (4)$ and $(1)-(2)= -6mn+4n-1 \ne 0$. Now, $(1)-(3)= 4mn(n-m)-11mn-2m^2+4m-4n+2 > 0$ for $m\le n-4,$ otherwise $(1)-(3) \le -n^2+9n-28 < 0$ for $m\ge n-3$. Thus $(1)\ne(3)$. Similarly, we can show that $(1) \ne (4), (5)$ and $(2) \ne (3), (4), (5)$. By Lemma~\ref{lem-G-e}, we have $\chi(H) \le \chi_{la}(H) \le 5$. Since  
 \end{proof}

\nt In~\cite[Theorems 3.9 \& 3.10]{LSN}, the authors proved that $\chi_{la}(C_m \vee K_n)=n+3$ for odd $m,n\ge 3$ and $\chi_{la}(C_m\vee K_n)=n+2$ for even $m\ge 4$, $n\ge 2$.

\begin{theorem}\label{thm-C2mVK2n-1} For $m\ge 2, n\ge 1$, $\chi_{la}(C_{2m}\vee K_{2n-1}) = 2n+1$. \end{theorem} 

\begin{proof}  Let $G=C_{2m} \vee K_{2n-1}$ be obtained from $C_{2m}\vee O_{2n-1}$ by adding $(n-1)(2n-1)$ edges joining $v_i,v_j$ for $1\le i < j \le 2n-1$. Suppose $n=1$, then $C_{2m}\vee K_1 = W_{2m}$, the wheel graph of order $2m+1$. In~\cite[Theorems 2.14 \& 5]{Arumugam, LNS}, the authors proved that $\chi_{la}(W_{2m}) = 3$ for $m\ge 2$. When $n=2$, $K_3 = C_3$ and the result follows from Theorem~\ref{thm-C2mVC2n-1}.

 \ms\nt We now consider $n\ge 3$. Let $g$ be the local antimagic labeling of $C_{2m}\vee O_{2n-1}$ as defined in the proof of Theorem~\ref{thm-C2mVO2n-1}. Let $h: E(K_{2n-1})\to [1,(n-1)(2n-1)]$ be a local antimagic labeling of $K_{2n-1}.$ Note that $h^+(v_j)$ are distinct for $1\leq j\leq 2n-1.$ Define $t: E(C_{2m}\vee K_{2n-1})\to [1,4mn+(n-1)(2n-1)] $ such that $t(e)=g(e)$ for $e\in E(C_{2m}\vee O_{2n-1})$ and $t(e)=h(e)+4mn$ for $e\in E(K_{2n-1})$. 

\ms\nt Thus, for $1\le i\le m$,  $t^+(u_{2i-1}) = g^+(u_{2i-1}) = m(4n^2-3n+3)+n$ and $t^+(u_{2i}) = g^+(u_{2i}) = m(4n^2 + 3n-1)+n+1$. Moreover, for $1\le j\le 2n-1$, $t^+(v_{j}) = h^+(v_j) + g^+(v_j) + 4mn(2n-2) = h^+(v_j) + m(4mn+2m+1)+ 4mn(2n-2)$. Since $h^+(v_j)$ are all distinct, we have $t^+(v_j)$ are all distinct too.

\ms\nt Clearly, $t^+(u_{2i-1}) < t^+(u_{2i})$ for $1\le i\le m$. Now, $t^+(v_j) - t^+(u_{2i-1}) = h^+(v_j) + m(4mn+2m+1)+ 4mn(2n-2)-m(4n^2-3n+3)-n = h^+(v_j) + m(4mn+2m-2+4n^2-5n)-n > 0$  for $m\ge 2, n\ge 3$.  Similarly, $t^+(v_j)\ne t^+(u_{2i})$. Thus, $t$ is a local antimagic labeling that induces $2n+1$ distinct vertex colors. Therefore, $\chi_{la}(C_{2m} \vee K_{2n-1}) \le 2n+1$. Since $\chi_{la}(C_{2m}\vee K_{2n-1}) \ge \chi(C_{2m}\vee K_{2n-1}) = 2n+1$, the theorem holds.
\end{proof}

\begin{theorem}\label{thm-K2nVC2m-1} For $m\ge 2, n\ge 1$, $\chi_{la}(K_{2n} \vee C_{2m-1})= 2n+3$. \end{theorem}

\begin{proof} Let $G=K_{2n} \vee C_{2m-1}$ such that $V(K_{2n}) = \{v_i\,|\,1\le i\le 2n\}$ and $V(C_{2m-1})=\{u_j: 1\le j \le 2m-1\}$ and $E(C_{2m-1})=\{u_ju_{j+1}: 1\le j \le 2m-1\}$ where $u_{2m}=u_1.$  Let $h$ be the local antimagic labeling that induces a $3$-coloring of $C_{2m-1}$ defined in the proof of~\cite[Theorem 2.8]{Arumugam}. Let $f$ be the local antimagic labeling of $K_{2n}$. Without loss of generality, we may assume $f^+(v_1) < f^+(v_3) < \cdots < f^+(v_{2n-1}) < f^+(v_2) < f^+(v_4) < \cdots < f^+(v_{2n})$. Define $g:E(G) \to [1,4mn+2m+n(2n-3)-1]$ such that $g(e) = h(e)$ if $e\in E(C_{2m-1})$, $g(v_iu_j) = x_{i,j} + (2m-1)$ for $1\le i\le 2n, 1\le j\le 2m-1$ where $x_{i,j}$ is the $(i,j)$-entry of a $(2n,2m-1)$-nearly magic rectangle that has $(2i-1)$-st row sum $\frac{(2m-1)(1+4mn-2n)-1}{2} = n(2m-1)^2+m-1$, $(2i)$-th row sum $\frac{(2m-1)(1+4mn-2n)+1}{2}=n(2m-1)^2+m$ and column sum $n(1+4mn-2n)$ (see~\cite{FSChai}). Finally, $g(e) = f(e) + (2n+1)(2m-1)$ if $e\in K_{2n}$.

\ms\nt Note that $h^+(u_1) = 3m-1$, $h^+(u_{2j-1}) = 2m-1$ for $j\ne 1$, and $h^+(u_{2j}) = 2m$. Hence,
\begin{enumerate}[(1)]
  \item $g^+(v_{2i-1}) = f^+(v_{2i-1}) + (2m-1)^2+ n(2m-1)^2 +m-1 $ for $1\le i\le n$, 
  \item $g^+(v_{2i}) = f^+(v_{2i})+(2m-1)^2+ n(2m-1)^2+m $ for $1\le i\le n$, 
  \item $g^+(u_1) = 3m-1 + 2n(2m-1) + 2n^2(2m-1)+n$,
  \item $g^+(u_{2j-1}) = 2m-1 + 2n(2n-1) + 2n^2(2m-1)+n$ for $2\le j\le m$,
  \item $g^+(u_{2j}) = 2m + 2n(2m-1) + 2n^2(2m-1)+n$ for $1\le j\le m-1$. 
\end{enumerate}   

\nt Clearly, $(2) > (1)$ and $(3) > (5) > (4)$. Moreover, $g^+(v_{2n}) > g^+(v_{2n-2}) > \cdots > g^+(v_2) > g^+(v_{2n-1}) > g^+(v_{2n-3}) > \cdots > g^+(v_1)$.  Since $f^+(v_1)\ge 1+2+\dots + (2n-1)+(2m-1)(2n+1)(2n-1)$,  we have $g^+(v_1)\ge n(2n-1)+(2m-1)(2n+1)(2n-1)+(2m-1)^2+ n(2m-1)^2+m-1$.  Thus $g^+(v_1) - g^+(u_1) = 2-8m+4m^2+n-8mn+4m^2n+4mn^2 = 4m(mn+m-2+n^2-2n)+n+2 > 0$. We can immediately conclude that $(1) > (3)$ for each possible $i$. Hence  $\chi_{la}(K_{2n} \vee C_{2m-1})\le 2n+3$ and since $\chi_{la}(G) \ge \chi(G)\ge 2n+3$, the theorem holds.
\end{proof}

\section{Conclusion and Open problems}

In this paper, we obtained the exact local antimagic chromatic number on the join of a path or a cycle with another graph. We end the paper with the following problems.

\begin{problem} For $m,n\ge 2$, determine $\chi_{la}(P_{2m-1}\vee O_n)$. \end{problem}

\begin{problem} For $m\ge 1, n\ge 2$, determine $\chi_{la}(P_{2m}\vee C_{2n})$. \end{problem} 

\begin{problem} For $m,n\ge 2$, determine $\chi_{la}(C_{2m}\vee C_{2n})$. \end{problem}

\begin{problem} For $m,n\ge 1$ and $m\ne n$, determine $\chi_{la}(C_{2m+1}\vee O_{2n})$. \end{problem} 

\begin{problem} For $m,n\ge 2$, determine $\chi_{la}((C_{2m}\vee C_{2n-1})-e)$ for $e\in E(C_{2n-1})$. \end{problem}

\begin{problem} For $m\ge 2, n\ge 1$, determine $\chi_{la}((C_{2m}\vee K_{2n-1})-e)$. \end{problem} 

\begin{problem} For $m\ge 2, n\ge 1$, determine $\chi_{la}(K_{2n} \vee C_{2m-1}-e)$. \end{problem}

\end{document}